\documentclass[11pt]{article}

\usepackage{mathpple}
\usepackage{latexsym}   % use all LaTeX fonts
\usepackage{amssymb}    % use all AMS fonts
\usepackage{amsmath}    % include some AMS-LaTeX functionality
\usepackage{amsthm}     % AMS style theorem and proof environments

\input{diagrams.tex}
\newarrow{to}---->
\newarrow{mto}|--->
\newarrow{l}-----
\newarrow{inject}C--->

\numberwithin{equation}{section}

\newtheorem{theorem}{Theorem}[section]
\newtheorem{lemma}[theorem]{Lemma}
\newtheorem{proposition}[theorem]{Proposition}
\newtheorem{corollary}[theorem]{Corollary}

\newenvironment{prf}[1]{\trivlist
\item[\hskip
\labelsep{\it #1.\hspace*{.3em}}]}{%~\hspace{\fill}~$\square$%
\endtrivlist}

\newtheorem{predefinition}[theorem]{Definition}
\newenvironment{definition}{\begin{predefinition}\rm}{\end{predefinition}}
\newtheorem{preremark}[theorem]{Remark}
\newenvironment{remark}{\begin{preremark}\rm}{\end{preremark}}
\newtheorem{prenotation}[theorem]{Notation}

\newtheorem{preexample}[theorem]{Example}

\newtheorem{preclaim}[theorem]{Claim}

\newtheorem{prequestion}[theorem]{Question}

 \makeatletter
\def\emppsubsection{\@startsection{subsection}{2}{\z@}{-3.25ex plus -1ex minus -.2ex}{-1em}{\bf}}

\makeatother

\newcommand \CM {{\cal M}}

\newcommand \ZZ {{\mathbb Z}}
\newcommand \NN {{\mathbb N}}

\newcommand  \FF {{\mathbb F}}
\newcommand \CC {{\mathbb C}}

\newcommand \dime {\mathop{\rm dim}}

\newcommand \Spec {\mathop{\rm Spec}}

\newcommand{\invlim}[1]{\lim_{\stackrel{\leftarrow}{#1}}}
\newcommand{\til}[1]{{\widetilde{#1}}}
\def\calm{\CM}

\def\smod{S_{{\rm mod}}}

\def\mmu{{\pmb \mu}}
\def\twiddle{{\sim}}
\def\th{{^{\rm th}}}
\newcommand{\st}[1]{\{#1\}}
\def\ra{\rightarrow}
\def\smooth{{\rm sm}}
\DeclareMathOperator{\sing}{Sing}
\DeclareMathOperator{\sym}{Sym}
\DeclareMathOperator{\spec}{Spec}
\DeclareMathOperator{\aut}{Aut}
\DeclareMathOperator{\gl}{GL}
\DeclareMathOperator{\pic}{Pic}
\DeclareMathOperator{\psp}{PSp}

\def\sp{{\mathop{\rm Sp}}}
\DeclareMathOperator{\gsp}{GSp}
\DeclareMathOperator{\Isom}{Isom}
\DeclareMathOperator{\End}{End}
\DeclareMathOperator{\diag}{diag}
\DeclareMathOperator{\Lie}{Lie}
\DeclareMathOperator{\mat}{Mat}
\def\std{{\mathop{\rm std}}}
\DeclareMathOperator{\sg}{SG}
\DeclareMathOperator{\su}{SU}
\DeclareMathOperator{\gu}{GU}

\DeclareMathOperator{\gal}{Gal}

\def\g{{\rm G}}
\def\sl{{\rm SL}}
\def\tensor{\otimes}
\def\integ{\mathbb Z}
\def\nat{\mathbb N}
\def\proj{\mathbb P}

\newcommand{\rest}[1]{|_{#1}}
\def\cx{\mathbb C}
\def\rat{\mathbb Q}
\def\ff{\mathbb F}
\def\real{\mathbb R}
\def\inject{\hookrightarrow}
\def\cross{\times}
\def\units{^\cross}
\DeclareMathOperator{\id}{id}
\def\iso{\cong}
\def\ang{\ip}
\def\mono{{\sf M}}
\newcommand{\ip}[1]{{{\langle #1 \rangle}}}

\def\calb{{\mathcal B}}
\def\calc{{\mathcal C}}
\def\calh{{\mathcal H}}
\def\call{{\mathcal L}}
\def\calo{{\mathcal O}}

\def\cals{{\mathcal S}}
\def\calsh{{\mathcal Sh}}
\def\calt{{\mathcal T}}

\def\dual{^\vee}
\def\inv{^{-1}}
\def\comp{\circ}
\renewcommand{\bar}[1]{{\overline{#1}}}

\newcommand{\bc}[1]{{\otimes #1}}

\newenvironment{alphabetize}{\begin{enumerate}

}{\end{enumerate}}
\newenvironment{romanize}{\begin{enumerate}

}{\end{enumerate}}

\title{The integral monodromy of hyperelliptic and trielliptic curves}
\author{Jeffrey D. Achter \& Rachel Pries
\footnote{The second author was partially supported by NSF grant DMS-04-00461.}}
\date{}

\begin{document}
\maketitle

\begin{abstract}
\noindent
We compute the $\integ/\ell$ and $\integ_\ell$ monodromy of every
irreducible component of the moduli spaces of hyperelliptic and trielliptic curves. 
In particular, we provide a proof that the $\integ/\ell$ monodromy of
the moduli space of hyperelliptic curves of genus 
$g$ is the symplectic group $\sp_{2g}(\integ/\ell)$.
We prove that the $\integ/\ell$ monodromy of the 
moduli space of trielliptic curves with signature $(r,s)$ is the
special unitary group
$\su_{(r,s)}(\integ/\ell\tensor\integ[\zeta_3])$.
{\bf keywords} monodromy, hyperelliptic, trigonal, moduli, Jacobian
\end{abstract}

\section{Introduction}

If $C \ra S$ is a relative smooth proper curve of genus $g \ge 1$ over an irreducible
base, then the $\ell$-torsion of the relative Jacobian of $C$  encodes
important information about the family.  Suppose $\ell$ is invertible
on $S$, and let $s \in S$ be a
geometric point.  The fundamental group $\pi_1(S,s)$ acts
linearly on the fiber $\pic^0(C)[\ell]_{s} \iso (\integ/\ell)^{2g}$, and
one can consider the mod-$\ell$ monodromy representation associated to $C$:
\begin{diagram}
\rho_{C \ra S, \ell}:\pi_1(S,s) & \rto & \aut(\pic^0(C)[\ell]_{s}) \iso
\gl_{2g}(\integ/\ell).
\end{diagram}
Let $\mono_\ell(C \ra S)$, or simply $\mono_\ell(S)$, be the image
of this representation. 
If a primitive $\ell\th$ root of unity is defined globally on $S$, then $\pic^0(C)[\ell]_{s}$ is equipped
with a skew-symmetric form $\ang{\cdot,\cdot}$ and $\mono_\ell(C \ra S) \subseteq
\sp(\pic^0(C)[\ell]_s,\ang{\cdot,\cdot}) \iso
\sp_{2g}(\integ/\ell)$. 
If $C \ra S$ is a sufficiently general family of curves, then
$\mono_\ell(C \ra S) \iso \sp_{2g}(\integ/\ell)$ \cite{delignemumford}.

In this paper, we compute $\mono_\ell(S)$ when $S$ is an
irreducible component of the moduli space of hyperelliptic or 
trielliptic curves and $C \ra S$ is the tautological curve.  
The first result implies that there is no restriction on the monodromy group
in the hyperelliptic case other than that it preserve the symplectic pairing.
As a trielliptic curve is a $\integ/3$-cover of a genus zero curve, 
the $\integ/3$-action
constrains the monodromy group to lie in a unitary group associated to $\integ[\zeta_3]$.
The second result implies that this is the only additional restriction in the trielliptic case.

\paragraph{Theorem \ref{thhe}}
{\it 
Let $\ell$ be an odd prime, and let $k$ be an algebraically closed field in which $2\ell$ is invertible.
For $g\ge 1$, $\mono_\ell(\calh_g\bc k)\iso
\sp_{2g}(\integ/\ell)$.}

\paragraph{Theorem \ref{thtri}}
{\it 
Let $\ell\ge 5$ be prime, and let $k$ be an algebraically closed field in which $3\ell$ is invertible.  Let
$\calt^{\bar\gamma}$ be any component of the moduli space of
trielliptic curves of genus $g\ge 3$. Then
$\mono_\ell(\calt^{\bar\gamma}\bc k) \iso
\sg_{(r_\gamma,s_\gamma)}(\integ/\ell)$ (where the latter is a unitary group defined
in \eqref{eqdefsg}).}

\medskip

We also prove that the $\ell$-adic monodromy group is 
$\sp_{2g}(\integ_\ell)$ in the situation of Theorem \ref{thhe} and is $\sg_{(r_\gamma,s_\gamma)}(\integ_\ell)$
in the situation of Theorem \ref{thtri}.

Theorem \ref{thhe} is an unpublished result of J.K. Yu and has already been used multiple times in the literature.
In \cite{chavdarov}, Chavdarov assumes this result to show that the numerator of the zeta function of
the typical hyperelliptic curve over a finite field is irreducible.
Kowalski also uses this result in a similar fashion \cite{kowalskisieve}.
The first author used Theorem \ref{thhe} to prove a conjecture of Friedman and
Washington on class groups of quadratic function fields \cite{achtercl}.

There are other results in the literature which are similar to Theorem \ref{thhe}
but which are not quite strong enough for the applications above.
A'Campo \cite[Th.\ 1]{acampo} computes the topological monodromy of $\calh_g \bc \CC$. 
On the arithmetic side, the $\rat_\ell$,
as opposed to $\integ_\ell$, monodromy of $\calh_g$
is computed in \cite[10.1.16]{katzsarnak}.  Combined with a theorem of
Larsen on compatible families of representations \cite[3.17]{larsenmax},
this shows that the mod-$\ell$ monodromy group 
of $\calh_g$ is maximal for a set of
primes $\ell$ of density one (as opposed to for all $\ell \ge 3$). 

There are results on $\rat_\ell$-monodromy of cyclic covers of the projective
line of arbitrary degree, e.g.,  \cite[Sec. 7.9]{katztwisted}.  Also,
in \cite[5.5]{fkv}, the authors prove that the projective representation
$\proj \rho_{C \ra S,\ell}$ is surjective for many
families of cyclic covers of the projective line.  
Due to a combinatorial hypothesis, their theorem does not apply to $\calh_g$
and applies to at most one component of the moduli space of
trielliptic curves for each genus, see Remark \ref{Rfkv}.   
See also work of Zarhin, e.g., \cite{zarhincyclic}.

As an application, for all $p \geq 5$, we show using \cite{chaioort01}
that there exist hyperelliptic and trielliptic curves
of every genus (and signature) defined over $\bar \ff_p$ whose Jacobians are absolutely simple.
In contrast with the applications above, 
these corollaries do not use the full strength of our results.
Related work can be found in \cite{HZhu} where the authors produce curves with absolutely simple 
Jacobians over $\ff_p$ under the restriction $g \leq 3$.

\paragraph{Corollary \ref{Chypabsirr}} 
{\it Let $p \not = 2$ and let $g\in\nat$.  Then there exists a
smooth hyperelliptic curve of genus $g$ defined over $\bar \ff_p$ whose Jacobian is
absolutely simple.}

\paragraph{Corollary \ref{Ctriabsirr}}
{\it Let $p \not = 3$.  Let $g \geq 3$ and let $(r,s)$ be a trielliptic signature for $g$
(Definition \ref{Dtrisig}).  
Then there exists a smooth trielliptic curve defined over $\bar \ff_p$ with genus $g$ and signature $(r,s)$
whose Jacobian is absolutely simple.}

\medskip 

Our proofs proceed by induction on the genus.
The base cases for the hyperelliptic family
rely on the fact that every curve of genus $g=1,2$ is hyperelliptic;
the claim on monodromy follows from the analogous assertion about the monodromy of $\calm_g$.
The base case $g=3$ for the trielliptic family involves a comparison with
a Shimura variety of PEL type, namely, the Picard modular variety.   
An important step is to show that the monodromy group does not change in the base cases when  
one adds a labeling of the ramification points to the moduli problem.

The inductive step is similar to the method used in \cite{ekedahlmono} 
and uses the fact that families of smooth hyperelliptic (trielliptic)
curves degenerate to trees of hyperelliptic (trielliptic) curves of lower genus.
The combinatorics of admissible degenerations require us 
to compute the monodromy exactly for the inductive step rather than up to isomorphism.

The inductive strategy using admissible degeneration developed here
should work for other families of curves, especially for more general
cyclic covers of the projective line.  The difficulty is in the direct
calculation of monodromy for the necessary base cases.

We thank C.-L.\ Chai, R.\ Hain, A.J.\ de Jong, E. Kani, and J. Kass.

\section{Moduli spaces of curves with $\integ/d$-action}

\subsection{Stable $\integ/d$-covers of a genus zero curve}

Let $G=\ZZ/d$ be a cyclic group of prime order $d$.  Let $G\units=G-\id_G$.
Let $S$ be an irreducible scheme over $\Spec \ZZ[1/d, \zeta_d]$.  Let
$k$ be an algebraically closed field equipped with a map
$\integ[1/d,\zeta_d] \ra k$.

Let $\psi:C \to S$ be a semi-stable curve.  In other words, $\psi$ is
flat and proper and the geometric fibers of $C$ are connected, reduced
curves whose only singularities are ordinary double points.  If $s \in
S$, let $C_s$ denote the fiber of $C$ over $s$.  Let $\sing_S(C)$
be the set of $z \in C$ for which $z$ is a singular point of the fiber
$C_{\psi(z)}$.

A {\it mark} $\Xi$ on $C/S$ is a closed subscheme of
$C-\sing_S(C)$ which is finite and \'etale over $S$.    The {\em
  degree} of $\Xi$ is the number of points in any geometric fiber of
$\Xi \ra S$.
A marked semi-stable curve $(C/S, \Xi)$ is {\it stably marked} if
every geometric fiber of $C$ satisfies the following condition: every
irreducible component of genus zero has at least three points which
are either in $\sing_S(C)$ or on the mark $\Xi$.

Consider a $G$-action $\iota_0: G \inject \aut_S(C)$ on $C$.  
Denote the ramification locus of $C \ra C/\iota_0(G)$ by $R$, and the
smooth ramification locus by $R_\smooth = R - (R\cap \sing_S(C))$.  We say
that $(C/S,\iota_0)$ is a {\em stable $G$-curve} if $C/S$ is a
semi-stable curve; if $\iota_0: G \inject \aut_S(C)$ is an action of $G$;
if $R_\smooth$ is a mark on $C/S$; and
if $(C/S,R_\smooth)$ is stably marked.  
We note that the definition implies that the {\it dihedral
nodes} of \cite[Def.\ 1.3]{ekedahlhurwitz} do not occur for $(C/S, \iota_0)$.

We say that a stable $G$-curve $(C/S, \iota_0)$ is {\it admissible} if
the following conditions are satisfied for every geometric point $z \in R\cap{\rm Sing}_S(C)$.  
Let $C_{z,1}$ and $C_{z,2}$ denote the two components of the formal
completion of $C_{\psi(z)}$ at $z$.
First, $\iota_0(1)$ stabilizes each branch $C_{z,i}$; 
%second, $z$ is a ramification point of the restriction of $\iota_0$ to $C_{z,i}$;
second, the characters of the action of $\iota_0$ on the tangent spaces
of $C_{z,1}$ and $C_{z,2}$ at $z$ are inverses. 

Throughout the paper, we suppose that $(C/S, \iota_0)$ is an admissible stable $G$-curve.  
We further assume throughout that $C/\iota_0(G)$ has arithmetic genus $0$.
Then $C/\iota_0(G)$ is also a stably marked curve \cite[Prop.\ 1.4]{ekedahlhurwitz}.  
The mark on $C/\iota_0(G)$ is the smooth branch locus
$B_\smooth$, which is the (reduced subscheme of) the image of $R_\smooth$ under the morphism $C
\ra C/\iota_0(G)$.
Let $r$ be the degree of $R_\smooth$.  By the Riemann-Hurwitz formula, the arithmetic
genus of each fiber of $C$ is $g=1-d+r(d-1)/2$.

Let $s$ be a geometric point of $S$ with residue field $k$
and let $a$ be a point of the fiber $R_{\smooth,s}$.  Then $G$ acts on the tangent space of $C_s$ at
$a$ via a character $\chi_a: G \ra k\units$.  In particular, there is
a unique choice of $\gamma_a \in (\integ/d)\units$ so that $\chi_a(1)
= \zeta_d^{\gamma_a}$.  
We say that $\gamma_a$ is the {\em canonical generator of inertia} at
$a$.  
The {\em inertia type} of $(C/S,\iota_0)$ is the multiset
$\st{\gamma_{a} : a \in R_{\smooth,s}}$.   It is independent of the choice of $s$.
By Riemann's existence theorem,
$\sum_{a \in R_{\smooth,s}} \gamma_a = 0 \in \integ/d$.

We say that a mark $\Xi$ has a labeling if $\Xi$ is an ordered disjoint
union of sections $S \ra C$. If $\Xi$ has degree $r$, we denote the
labeling by $\eta:\st{1, \ldots, r} \ra \Xi$.
%A {\em labeling} of $\Xi$ is a choice $\eta$ of $r$ ordered disjoint sections such
%that $\Xi$ is the disjoint union of the image of the sections. 
A labeling of an admissible stable $G$-curve $(C/S,\iota_0)$ is a
labeling $\eta$ of $R_\smooth$.
There is an induced labeling $\eta_0:\{1, \ldots, r\} \ra B_\smooth$.

If $(C/S,\iota_0,\eta)$ is a labeled $G$-curve, the {\em class vector}
is the map of sets $\gamma:\st{1, \ldots, r} \ra G\units$ such that
$\gamma(i) = \gamma_{\eta(i)}$.  We frequently write $\gamma = (\gamma(1),
\ldots, \gamma(r))$.  If $\gamma$ is a class vector, we denote
its inertia type by $\bar\gamma: G\units \ra \integ_{\ge 0}$ 
where $\bar\gamma(h) = \#\gamma\inv(h)$ for all $h\in G$.

\subsection{Moduli spaces}

We define moduli functors on the category of schemes
over $\spec\integ[1/d, \zeta_d]$ by describing their $S$-points:
\begin{description}
\item{$\overline{\calm}_G$} parametrizes admissible stable $G$-curves $(C/S, \iota_0)$. 
\item{$\til\calm_G$} parametrizes labeled admissible stable $G$-curves $(C/S, \iota_0, \eta)$. 
\item{$\til\calm_{g,r}$} parametrizes triples $(C/S, \Xi,\eta)$ where
  $C/S$ is a semi-stable curve of genus $g$, $\Xi$ is a mark of degree
  $r$ on $C$ such that $(C/S,\Xi)$ is stably marked, and $\eta$ is a
  labeling of $\Xi$.
\end{description}
Each functor is represented by an algebraic stack, and we use the
same letter to denote both a moduli functor and its representing
stack.
For each of these moduli spaces $\calm$, we let $\calm^{\circ}$ denote the open substack 
whose objects parametrize smooth curves of the appropriate type.
To work with fibers
of the structural map $\calm \ra \spec
\integ[1/d,\zeta_d]$ we 
write $\calm\bc k$ for $\calm\cross_{\spec \integ[1/d,\zeta_d]}
\spec k$.

\begin{lemma} \label{Lmoduli}
The moduli spaces $\til\calm_G$ and $\bar\calm_G$ are smooth, proper
Deligne-Mumford stacks over $\spec\integ[1/d,\zeta_d]$.  The subspaces
$\til\calm_G^\circ$ and $\bar{\calm}^\circ_G$ are open and dense in
$\til\calm_G$ and $\bar\calm_G$, respectively.  
\end{lemma}

\begin{proof}
The moduli spaces $\til\calm_G$
and $\bar\calm_G$ are algebraic stacks.  Since the
automorphism group scheme of a (labeled) admissible stable $G$-curve
over an algebraically closed field is \'etale, $\bar\calm_G$ and
$\til\calm_G$ are Deligne-Mumford stacks \cite[8.1]{lmbstacks}.

The local deformation problem
for a semi-stable curve with $G$-action is formally smooth
\cite[2.1]{ekedahlhurwitz}, so $\bar\calm_G$ is smooth.  By
\cite[2.3]{ekedahlhurwitz}, the degree of the smooth
ramification divisor of a stable admissible $G$-curve is locally constant, so
$\til\calm_G$ is smooth, too.

We use the valuative criterion to show that these stacks are proper \cite[7.3]{lmbstacks}.
Let $\calo_K$ be a
complete discrete valuation ring with field of fractions $K$, and suppose
$(C,\iota_0) \in \bar\calm_G(K)$.  By the stable reduction theorem,
there exists a finite extension $K'/K$ so that the curve $C$ extends
as a stable curve to $\calo_{K'}$, as does the admissible $G$-action $\iota_0$.  Possibly after a
further extension of the base, one may blow up the special fiber of
$C/\calo_{K'}$ to remove any dihedral nodes
\cite[p.195]{ekedahlhurwitz}, so that no singular point of $C$ is in
the closure of $R_\smooth$ and the resulting curve is an
admissible stable $G$-curve.  This  shows that $\bar\calm_G$ is
proper.  After a finite base change, a labeling of $R_\smooth$
also extends uniquely.  Thus, 
$\til\calm_G$ is also proper.

The openness and density of
$\til\calm_G^\circ$ in $\til\calm_G$ and of $\bar\calm_G^\circ$ in
$\bar\calm_G$ follow from the fact that an admissible stable
$G$-curve is equivariantly smoothable \cite[2.2]{ekedahlhurwitz}.
\end{proof}

Consequently, every connected component of $\bar\calm_G$ or $\til\calm_G$ is irreducible.

Let $\gamma: \st{1, \ldots, r} \ra (\integ/d)\units$ be a class vector
of length $r=r(\gamma)$ for $G$.  Let $g(\gamma) :=1 - d+r(\gamma)(d-1)/2$.
Let $\til\calm_G^\gamma$ be the substack of $\til\calm_G$ 
for which $(C/S,\iota_0,\eta)$ has class vector $\gamma$. 
Let $\bar\calm_G^{\bar\gamma}$ be the substack of $\bar\calm_G$ for
which $(C/S,\iota_0)$ has inertia type $\bar\gamma$.

\begin{lemma} \label{Lirreducible}
The moduli space $\til\calm_G^\gamma$ is irreducible. 
\end{lemma}

\begin{proof}
Since $\til\calm_G$ is proper and smooth over $\spec \integ[1/d,
\zeta_d]$, it is sufficient to prove that $\til\calm_{G}^\gamma \bc \cx$
is irreducible \cite[IV.5.10]{faltingschai}.  
Since $\til\calm_G^{\gamma,\circ}\bc \cx$ is open and dense in $\til
\calm_{G}^\gamma\bc \cx$, it suffices to prove that
$\til\calm_{G}^{\gamma,\circ}\bc \cx$ is irreducible.
Consider the functor $\beta_\gamma: \til \calm_G^\gamma \ra \til\calm_{0,r(\gamma)}$. 
On $S$-points, this functor takes the isomorphism class of
$(C,\iota_0, \eta)$ to the isomorphism class of $(C/\iota_0(G),
B_\smooth,\eta_0)$, where $\eta_0$ is
the induced labeling of $B_\smooth$.
By Riemann's existence theorem \cite[Section 2.2]{Vbook},
$\beta^{\circ}_{\gamma}\bc \cx: \til \calm^{\gamma, \circ}_G\bc\cx \ra \til\calm^{\circ}_{0,r}\bc\cx$ is an isomorphism.
The statement follows since $\til\calm^{\circ}_{0,r}\bc\cx$ is irreducible.
\end{proof}

If two class vectors $\gamma$ and $\gamma'$ yield the same inertia type, so that
$\bar{\gamma} = \bar{\gamma}'$, then there is a permutation $\varpi$ of
$\st{1, \ldots, r}$ such that $\gamma' = \gamma \comp \varpi$.  This
relabeling of the branch locus yields an isomorphism
\begin{diagram}[LaTeXeqno] \label{Erelabel}
\label{diagrelabellocus}
\til\calm^{\gamma}_G & \rto^\twiddle & \til \calm^{\gamma \comp \varpi}_G.
\end{diagram}

Suppose $\gamma$ and $\gamma'$ differ by an
automorphism of $G$, so that there exists $\sigma\in \aut(G)$ such
that $\gamma'= \sigma \comp \gamma$.  This relabeling of
the $G$-action yields an isomorphism 
\begin{diagram}[LaTeXeqno]\label{Eaut}
\label{diagrelabelG}
\til\calm^{\gamma} &\rto^\twiddle & \til \calm^{\sigma\comp \gamma}_G.
\end{diagram}

\begin{lemma}\label{lemforgetgalois}
The forgetful functor $\til\calm_G^\gamma \ra
\bar\calm_G^{\bar\gamma}$ is \'etale and Galois.
\end{lemma}

\begin{proof}
The map is \'etale since any admissible stable
$G$-curve admits a labeling of $R_\smooth$ \'etale-locally on the base.
Moreover, the set of labelings of a
fixed stable $G$-curve with rational smooth ramification locus is a
torsor under the subgroup of $\sym(r)$ consisting
of those $\varpi$ for which $\gamma = \gamma\comp \varpi$.  Therefore,
$\til\calm_G^\gamma \ra \bar\calm_G^{\bar\gamma}$ is Galois.
\end{proof}

\subsection{Degeneration}

For $i=1,2,3$, let $\gamma_i$ denote a class vector with length $r_i$
and let $g_i=g(\gamma_i)$.
There is a clutching map 
\begin{diagram}
\label{diagrelabelG2}
\til\calm_{g_1,r_1} \cross \til \calm_{g_2,r_2} & \rto &
\til\calm_{g_1+g_2, r_1+r_2-2}
\end{diagram}
which is a closed immersion \cite[3.9]{knudsen2}.  On $S$-points, this map
corresponds to gluing $C_1$ and $C_2$ together over $S$ by
identifying the last section of $C_1$ with the first section of $C_2$.

There is a forgetful functor $\til \calm_G^\gamma \ra \til \calm_{g(\gamma), r(\gamma)}$ taking the 
isomorphism class of $(C/S, \iota_0, \eta)$ to the isomorphism class of $(C/S, R_\smooth,\eta)$.
This functor is finite-to-one since ${\rm Aut}_S(C, \iota_0)$ is finite \cite[1.11]{delignemumford}.
The composition 
\begin{diagram}
\til \calm^{\gamma_1}_G
\cross \til \calm^{\gamma_2}_G &\rto &  \til \calm_{g_1,r_1} \cross \til
\calm_{g_2, r_2} & \rto& \til \calm_{g_1+g_2, r_1+r_2-2}
\end{diagram}
allows us to glue two labeled $G$-curves $(C_i/S, \iota_{0,i}, \eta_i)$ together to obtain a labeled 
$G$-curve $C/S$ with genus $g_1+g_2$ and class vector
$$\gamma=(\gamma_1(1), \ldots, \gamma_1(r_1-1), \gamma_2(2), \ldots, \gamma_2(r_2)).$$  
Moreover, $C/S$ is
equivariantly smoothable if and only if the $G$-action is admissible,
i.e., if and only if $\gamma_1(r)$ and $\gamma_2(1)$ are inverses \cite[2.2]{ekedahlhurwitz}.  
In this situation, we say that
{\it $(\gamma_1,\gamma_2)$ deforms to $\gamma$} or that {\it $\gamma$
degenerates to $(\gamma_1,\gamma_2)$}, and write 
\begin{diagram}
\til\calm^{\gamma_1}_G \times \til\calm^{\gamma_2}_G & \rto &\til\calm^{\gamma}_G.
\end{diagram} 

The clutching maps admit a generalization to maps 
\begin{diagram}
\til \calm_{g_1,r_1} \cross \til \calm_{g_2,r_2} \cross \til
\calm_{g_3,r_3} &\rto & \til \calm_{g_1+g_2+g_3, r_1+r_2 +r_3- 4}.
\end{diagram}

Assume that $(\gamma_1, \gamma_2)$ deforms to $\gamma_L$, 
that $(\gamma_2, \gamma_3)$ deforms to $\gamma_R$, 
and that both $(\gamma_L, \gamma_3)$ and $(\gamma_1, \gamma_R)$ deform to $\gamma$.  
Then we have a commutative diagram of maps:
\begin{diagram}[LaTeXeqno]
\label{diagclutch}
\til\calm_G^{\gamma_1} \cross \til\calm_G^{\gamma_2} \cross \til\calm_G^{\gamma_3} & 
\rto & \til\calm_G^{\gamma_1} \cross \til\calm_G^{\gamma_R} \\
\dto & \rdto & \dto \\
\til\calm_G^{\gamma_L} \cross \til \calm_G^{\gamma_3} & \rto &
\til\calm_G^{\gamma}.
\end{diagram}

In the situation above, suppose $g_1=g_3=1$.  The image $\Delta_{1,1}$ of 
$\til\calm_G^{\gamma_1} \cross \til\calm_G^{\gamma_2} \cross \til\calm_G^{\gamma_3}$
in $\til \calm_G^\gamma$
is part of the boundary of $\til\calm_G^\gamma$.  

We say that a class vector $\gamma$ {\it degenerates to $\Delta_{1,1}$} 
if there are $\gamma_1$, $\gamma_2$ and $\gamma_3$ as in
\eqref{diagclutch} so that $g_1=g_3=1$.   
Furthermore, we say that an inertia type $\overline{\gamma}$ {\it degenerates to $\Delta_{1,1}$} 
if there exists some class vector $\gamma'$ with inertia type $\overline{\gamma}' =\overline{\gamma}$ 
so that $\gamma'$ degenerates to $\Delta_{1,1}$.
Using Equations \eqref{diagrelabellocus} and
\eqref{diagrelabelG}, if $\overline{\gamma}$ degenerates to $\Delta_{1,1}$,
then every component of $\til \calm_G$ lying over
$\bar\calm_G^{\overline{\gamma}}$ degenerates to $\Delta_{1,1}$
(although this may require identifying different ramification points
of $C_1$, $C_2$  and $C_3$ in the clutching maps).

\subsection{The case of hyperelliptic curves: $d=2$} \label{subsechecomp}

Let $C/S$ be a $\integ/2$-curve of genus $g$.  Then $C/S$ is
hyperelliptic, and $\iota_0(1)$ is the hyperelliptic involution. Over
every geometric point $s$ of
$S$, $C_s \ra C_s/\iota_0(\integ/2)$ is ramified at $2g+2$ smooth points.

There is a unique class vector $(1, \ldots, 1)$ for $\integ/2$ of
length $2g+2$ for each $g\in\nat$.  By Lemma \ref{Lirreducible}, there
is a unique component of $\til\calm_{\integ/2}$
parametrizing labeled hyperelliptic curves of genus $g$.  We denote
this component by $\til\calh_g$.  Similarly, $\bar\calm_{\integ/2}$
has a unique component $\bar\calh_g$ which parametrizes
hyperelliptic curves of genus $g$.  Let $\calh_g=\bar\calh_g^\circ$ denote
the moduli space of smooth hyperelliptic curves of genus $g$.

\begin{lemma}
\label{lemdegenhe}
If $g \geq 3$, then $\til{\calh}_g$ degenerates to $\Delta_{1,1}$.
\end{lemma}

\begin{proof}
The class vector for $\til\calh_g$ is 
$\gamma=(1, \ldots, 1)$ with length $2g+2$.
Let $\gamma_1=\gamma_3=(1,1,1,1)$ and let $\gamma_2=(1, \ldots, 1)$ with length $2g-2$.
\end{proof}

\begin{remark}
When the boundary components $\Delta_1$ and $\Delta_2$ are distinct,
the inductive argument in Theorem \ref{thhe} can be revised 
to rely on monodromy groups only up to 
isomorphism.  Unfortunately, when $g = 3$ one has $\Delta_1 =
\Delta_2$, so that we must calculate monodromy groups exactly.
\end{remark}

\subsection{The case of trielliptic curves: $d=3$} \label{subsectricomp}

A {\it trielliptic} curve is a $\ZZ/3$-curve $(C/S, \iota_0)$ so that
$C/\iota_0(\ZZ/3)$ has genus zero.  
The curve $C/S$ is sometimes called a cyclic trigonal curve.
Over every geometric point
$s$ of $S$, the $\integ/3$-cover $C_s \ra C_s / \iota_0(\integ/3)$ is
ramified at $g+2$ smooth points.

\subsubsection{Components of the trielliptic locus}

The moduli space of trielliptic curves $(C/S, \iota_0)$ with $C$ of genus $g$ 
is not connected even for fixed $g$. 
This reflects the different possibilities for the inertia type, 
or equivalently for the signature, which we now describe.

By Lemma \ref{Lirreducible}, the components of $\til
\calm_{\integ/3}$ parametrizing labeled trielliptic curves of genus $g$ are in
bijection with the maps $\gamma: \st{1, \ldots, g+2} \ra
(\integ/3)\units$ such that $\sum \gamma(i) = 0 \in \integ/3$. We
denote these components by $\til\calt^\gamma$, and their projections
to $\bar\calm_{\integ/3}$ by $\bar\calt^{\overline{\gamma}}$.  Let 
$\calt^{\bar\gamma}={\bar\calt^{\bar\gamma\circ}}$ denote the
moduli space of smooth trielliptic curves with inertia type $\bar\gamma$.

\begin{lemma} \label{Lgamma}
The set of possible inertia types $\overline{\gamma}$ for a trielliptic curve $(C/S, \iota_0)$ of genus $g$
is in bijection with the set of pairs of integers $(d_1,d_2)$ so that 
$d_1,d_2 \geq 0$, $d_1+d_2=g+2$, 
and $d_1+2d_2 \equiv 0 \bmod 3$.
\end{lemma}

\begin{proof}
Let $d_1=\bar \gamma(1)$ and $d_2=\bar \gamma(2)$.  The claim follows from earlier remarks.
\end{proof}

\begin{remark} \label{Rfkv}
The result \cite[5.5]{fkv} applies to
trielliptic curves only when $d_1d_2 = 0$.  These inertia types occur
for curves of genus
  $3g_1+1$ having an equation of the form $y^3=f(x)$ for some
  separable polynomial $f(x)$ of degree $3g_1+2$.
\end{remark}

Let $S$ be an irreducible scheme over $\spec\integ[1/3, \zeta_3]$.  Then
$\calo_S\tensor \integ[1/3,\zeta_3] \iso \calo_S \oplus \calo_S$.  We
choose the isomorphism so that the first component has the
given structure of $\calo_S$ as a $\integ[1/3,\zeta_3]$-module.

Consider a trielliptic curve $(\psi: C \ra S, \iota_0)$.
The sheaf of relative
one-forms $\psi_*\Omega^1_{C/S}$ is a locally free
$\integ[1/3,\zeta_3]\tensor \calo_S$-module of some rank $(r,s)$, where $r+s = g$.  We
call $(r,s)$ the {\em signature} of $(C/S, \iota_0)$.  The signature is locally
constant on $S$, so we may calculate it at any geometric point of $S$.
If $(C/k,\iota_0)$ is a trielliptic curve, 
then $H^0(C,\Omega^1_C)$ decomposes as a direct sum $W_1 \oplus W_2$ where
$\omega \in W_j$ if $\zeta_3 \circ \omega = \zeta_3^j \omega$.
The signature of $(C/k, \iota_0)$ is $(r,s) = (\dime(W_1),\dime(W_2))$.

\begin{lemma} \label{Lsig}
The signature and inertia type of a trielliptic curve of genus $g$ are related as follows:
$\overline{\gamma}(1)=2r-s+1$ and $\overline{\gamma}(2)=2s-r+1$.
There exists a trielliptic curve $(C/k, \iota_0)$ of genus $g$ with signature $(r,s)$ if and only if 
$r,s \in \ZZ$, $r+s=g$, and $(g-1)/3 \leq r,s \leq (2g+1)/3$.
\end{lemma}

\begin{proof}
Let $(C/k, \iota_0)$ be a trielliptic curve with inertia type $\overline{\gamma}$.
For simplicity, let $d_1=\overline{\gamma}(1)$ and let $d_2=\overline{\gamma}(2)$ 
and let $N=(d_1+2d_2)/3$.  There is an equation for $C$ of the form
$$y^3=\prod_{i=1}^{d_1}(x-a_i)\prod_{j=1}^{d_2}(x-b_j)^2.$$

Consider the differential 
$\omega=\prod_{i=1}^{d_1}(x-a_i)^{n_i}\prod_{j=1}^{d_2}(x-b_j)^{n_j'}dx/y^m$ 
with $n_i, n_j' \geq 0$ and $m \geq 1$.
By \cite[Thm.\ 3]{koo91}, $\omega$ is holomorphic if and only if 
$3n_i \geq m-2$, $3n_j' \geq 2m-2$, and $mN \geq \sum_{i=1}^{d_1} n_i +\sum_{j=1}^{d_2} n_j' +2$. 
Let $g(x)=\prod_{j=1}^{d_2}(x-b_j)$.
Thus we have the following set of linearly independent holomorphic differentials:
$$\{dx/y, \ldots, (x-a_1)^{N-2}dx/y, g(x)dx/y^2, \ldots, (x-a_1)^{2N-2-d_2}g(x)dx/y^2\}.$$
This set is a basis of $H^0(C,\Omega^1_C)$ since its cardinality is $3N-d_2-2=g$.
Also, for any $h(x) \in k[x]$, $\zeta_3 \circ h(x)dx/y^m = \zeta_3^{-m} h(x) dx/y^m$.
It follows that the signature of $(C/k, \iota_0)$ is $(r,s)=(g-N+1,N-1)$.

It follows that $d_1=2r-s+1$ and $d_2=2s-r+1$.     
Now $r,s \in \ZZ$ if and only if $d_1+2d_2 \equiv 0 \bmod 3$.
Also, $d_1+d_2=g+2$ if and only if $r+s=g$.
The conditions $d_1,d_2 \geq 0$ and $(g-1)/3 \leq r,s \leq (2g+1)/3$ are equivalent.
The second claim then follows from Lemma \ref{Lgamma}.
\end{proof}

\begin{definition} \label{Dtrisig}
Let $g \in \NN$.  A {\it trielliptic signature for $g$} is a pair $(r,s)$ with 
$r,s \in \ZZ$, $r+s=g$, and $(g-1)/3 \leq r,s \leq (2g+1)/3$.
\end{definition}

As in \eqref{diagrelabelG}, if $(C/S, \iota_0)$ is a trielliptic curve then so is $(C/S, \iota'_0)$  
where $\iota'_0(1)=\iota_0(2)$. 
Replacing $\iota_0$ with $\iota'_0$ exchanges the values of $d_1$ and $d_2$ 
and the values of $r$ and $s$.

\subsubsection{Degeneration of the trielliptic locus}

We show that every  component $\til\calt^\gamma$ with $g(\gamma) \geq 4$
degenerates to $\Delta_{1,1}$.   
 
Recall that there is a unique elliptic curve $E$ which admits a $\ZZ/3$-action $\iota_0$.  
The class vector of $(E, \iota_0)$ is either $(1,1,1)$ or $(2,2,2)$.

\begin{proposition} \label{Pdeg311}\label{propdegentri}
If $g(\gamma) \geq 4$, then $\til \calt^\gamma$ degenerates to $\Delta_{1,1}$.
\end{proposition}

\begin{proof}
Let $g=g(\gamma) \geq 4$ and 
let $\gamma:\{1, \ldots,g+2\} \ra (\ZZ/3)\units$ be a class vector for $G=\ZZ/3$.
It suffices to show that there exist class vectors $\gamma_i$ for $i=1,2,3$
with $g(\gamma_1)=g(\gamma_3)=1$ as in \eqref{diagclutch}.
In particular, we require that $\gamma_1$ and $\gamma_3$ are either $(1,1,1)$ or $(2,2,2)$.

By Equations \eqref{Erelabel} and \eqref{Eaut}, we can reorder the
values $\gamma(i)$ or replace each $\gamma(i)$ with $-\gamma(i)$.  
Thus it is sufficient to restrict attention to the inertia type $\overline{\gamma}$,
or to the corresponding signature $(r,s)$ by Lemma \ref{Lsig}.

Suppose $\gamma_1=\gamma_3=(1,1,1)$.  Then we need $\gamma_2:\{1,
\ldots, g\} \to (\ZZ/3) \units$ with $\gamma_2(1)=\gamma_2(g)=2$ and
$\overline{\gamma}_2(1)=\overline{\gamma}(1)-4$ and
$\overline{\gamma}_2(2)=\overline{\gamma}(2)+2$.  In other words, we
need $\gamma_2$ to have four fewer points with inertia generator 1 and
two more points with inertia generator 2 than $\gamma$ does.  This can
be achieved by decreasing $r$ by $2$.  Similarly, the case
$\gamma_1=\gamma_3=(2,2,2)$ can be achieved by decreasing $s$ by $2$.

Suppose $\gamma_1=(1,1,1)$ and $\gamma_2=(2,2,2)$.
Then we need $\gamma_2:\{1, \ldots, g\} \to (\ZZ/3) \units$ with $\gamma_2(1)=2$ and $\gamma_2(g)=1$ and 
$\overline{\gamma}_2(1)=\overline{\gamma}(1)-1$ and $\overline{\gamma}_2(2)=\overline{\gamma}(2)-1$.
In other words, we need $\gamma_2$ to have one fewer point with inertia generator 1 and 
one fewer point with inertia generator 2 than $\gamma$ does.  
This can be done by decreasing both $r$ and $s$ by $1$.

The next table shows that, for each inertia type $\overline{\gamma}$, there is a choice of
$\gamma_1, \gamma_2, \gamma_3$ satisfying the numerical constraints. 

\[\begin{tabular}{|l|l|l|l|l|}
\hline
Component $\gamma$ & degenerates so that & and $\gamma_1$ is & and $\gamma_3$ is & under this \\
with signature & $\gamma_2$ has signature & &  & condition\\
\hline
\hline
$(r,s)$ & $(r-2,s)$& $(1,1,1)$ & $(1,1,1)$ & $r \geq 2$ \\
\hline
$(r,s)$ & $(r,s-2)$ & $(2,2,2)$ & $(2,2,2)$ & $s \geq 2$ \\
\hline
$(r,s)$ & $(r-1,s-1)$ & $(1,1,1)$ & $(2,2,2)$ & $r \geq 1$, $s \geq 1$ \\
\hline
\end{tabular}
\]
\end{proof}

\begin{remark}
  When $\til \calt^\gamma$ degenerates to both $\Delta_1$ and
  $\Delta_2$, the inductive argument in Theorem \ref{thtri} can be
  revised to rely on monodromy groups only up to isomorphism.
  Unfortunately, $\til{\calt}^\gamma$ does not degenerate to $\Delta_2$
  when $\bar\gamma(1)\bar\gamma(2)=0$.    For $g(\gamma) \geq
  4$, all other components $\til{\calt}^{\gamma}$ degenerate to
  $\Delta_2$.

\end{remark}

\section{Monodromy groups}

\subsection{Definition of monodromy}

Let $(X/S,\phi)$ be a principally polarized abelian scheme of
relative dimension $g$ over an irreducible base.  If $\ell$ is a
rational prime invertible on $S$, then the $\ell$-torsion $X[\ell]$ of
$\ell$ is an \'etale cover of $S$ with geometric fiber isomorphic to
$(\integ/\ell)^{2g}$.  
Let $s$ be a geometric point of $S$.  The fundamental group $\pi_1(S,s)$ acts
linearly on the $\ell$-torsion of $X$.
This yields a representation
\begin{diagram}
\rho_{X \ra S, s,\ell}: \pi_1(S,s) & \rto & \aut(X[\ell]_s) \iso \gl_{2g}(\integ/\ell).
\end{diagram}
The cover $X[\ell] \ra S$ both determines and is determined by the representation 
$\rho_{X \ra S, s,\ell}$.  The image of $\rho_{X \ra S, s,
  \ell}$ is the {\it mod-$\ell$ monodromy} of $X \ra S$ and we denote it by
$\mono_\ell(X \ra S, s)$, or by $\mono_\ell(S,s)$ if the choice of
abelian scheme is clear.  The isomorphism class of the
$\mono_\ell(S,s)$ is independent of the choice of base point $s$, and we
denote it by $\mono_\ell(S)$.

Let $X\dual$ be the dual abelian scheme.  There
is a canonical pairing $X[\ell] \cross X\dual[\ell] \ra
\mmu_{\ell,S}$, where $\mmu_{\ell,S} := \mmu_\ell \cross S$ is
the group scheme of $\ell\th$ roots of unity.
The polarization $\phi$ induces an isomorphism $X \ra X\dual$, and
thus a skew-symmetric pairing $X[\ell] \cross X[\ell] \ra \mmu_{\ell,S}$.
Because the polarization is defined globally, the image of monodromy
$\mono_\ell(X \ra S, s)$ is contained in the group of symplectic
similitudes of $(X[\ell]_s,
\ang{\cdot,\cdot}_\phi)$, which is isomorphic to
$\gsp_{2g}(\integ/\ell)$.  Moreover, if a primitive $\ell^{{\rm th}}$ root of
unity exists globally on $S$, then $\pi_1(S,s)$ acts trivially on
$\mmu_{\ell,S}$ and $\mono_\ell(X \ra S,s) \subseteq
\sp(X[\ell]_s,\ang{\cdot,\cdot}_\phi) \iso \sp_{2g}(\integ/\ell)$.

Similarly, the cover $X[\ell^n] \ra S$ defines a monodromy representation 
with values in $\aut(X[\ell^n]_s) \iso\gl_{2g}(\integ/\ell^n)$. Taking
the inverse limit over all $n$, we obtain a continuous
representation on the Tate module of $X$, 
\begin{diagram}
\rho_{X \ra S, \integ_\ell, s}: \pi_1(S,s) & \rto & \invlim n
\aut(X[\ell^n]_s) \iso \gl_{2g}(\integ_\ell).
\end{diagram}
We denote the image of this representation by $\mono_{\integ_\ell}(X
\ra S, s)$, and its isomorphism class by $\mono_{\integ_\ell}(X \ra
S)$ or $\mono_{\integ_\ell}(S)$.  Again, there is an inclusion 
$\mono_{\integ_\ell}(X \ra S) \subseteq \gsp_{2g}(\integ_\ell)$.  If
$F$ is a field,  let $F_{\ell^\infty} = F(\mmu_{\ell^\infty}(\bar F))$.
If
$S$ is an $F$-scheme, then $\mono_{\integ_\ell}(X \ra S,
s)/ \mono_{\integ_\ell}(X\bc{\bar F} \ra S \bc{\bar F}, s) \iso
\gal(F_{\ell^\infty}/F)$.  Finally, let $\mono_{\rat_\ell}(X\ra S, s)$
be the Zariski closure of $\mono_{\integ_\ell}(X \ra S, s)$ in
$\gl_{2g}(\rat_\ell)$.

Now suppose that $\psi:C \ra S$ is a relative proper semi-stable curve.
Let $\pic^0(C) := \pic^0_{C/S}$ be the neutral component of the
relative Picard functor of $C$ over $S$.  Since $C/S$ is semi-stable,
$\pic^0(C)$ is a semiabelian scheme \cite[9.4.1]{blr}.  Suppose that
there is at least one geometric point $s$ such that the fiber
$\pic^0(C_s)$ is an abelian variety.  (This is true if some $C_s$ is a tree
of smooth curves.)  Then there is a nonempty open
subscheme $S^*$ of $S$ such that $\pic^0(C\rest{S^*})$ is an abelian scheme
over $S^*$.  We define the mod-$\ell$ and $\integ_\ell$ monodromy
representations of $C$ to be those of $\pic^0(C\rest{S^*}) \ra S^*$.
(Alternatively, these may be constructed as the restrictions of
$R^1\psi_*\mmu_{\ell,S}$ and $R^1\psi_*\mmu_{\ell^\infty,S}$
to the largest subscheme of $S$ on which these sheaves are unramified.)
Thus, $\mono_\ell(C \ra S, s) = \mono_\ell(\pic^0(C\rest{S^*}) \ra S^*, s)$,
and we denote this again by $\mono_\ell(S,s)$ if the curve is clear
and by $\mono_\ell(S)$ if the base point is suppressed.

The moduli spaces $\bar\calm_G$ and $\til\calm_G$ are Deligne-Mumford
stacks, and we employ a similar formalism for \'etale covers of stacks
\cite{noohi}.  Let $\cals$ be a connected Deligne-Mumford stack.  The
category of Galois \'etale covers of $\cals$ is a Galois category in
the sense of Grothendieck, and thus there is an \'etale fundamental
group of $\cals$.  More precisely, let $s\in \cals$ be a geometric
point.  Then there is a group $\pi_1(\cals,s)$ and an equivalence of
categories between finite $\pi_1(\cals,s)$-sets and finite \'etale
Galois covers of $\cals$.  If $\cals$ has a coarse moduli space
$\smod$, then $\pi_1(\cals,s)$ is the extension of $\pi_1(\smod,s)$ by
a group which encodes the extra automorphism structure on the moduli
space $\smod$ \cite[7.11]{noohi}. If $X \ra \cals$ is a family of
abelian varieties, we again let $\mono_\ell(X\ra \cals,s)$ be the
image of $\pi_1(\cals, s)$ in $\aut(X[\ell]_s)$.  

Let
$\calc^\gamma$ be the tautological labeled curve over
$\til\calm^\gamma_G$.  By the mod-$\ell$ or $\integ_\ell$ monodromy of
$\til\calm_G^\gamma$ we mean that of $C^\gamma \ra
\til\calm_G^\gamma$.

\subsection{Degeneration and monodromy} \label{Sinduct}

In this section, we deduce information on monodromy groups from the 
degeneration of moduli spaces to the boundary.
We will use this to prove inductively that the mod-$\ell$
monodromy groups of each $\til\calh_g$ and each $\til\calt^\gamma$
are as large as possible.  

\begin{lemma}
\label{lemclutchmono}
Let $k$ be an algebraically closed field in which $d\ell$ is invertible.
Suppose that the pair $(\gamma_1,\gamma_2)$ deforms to $\gamma$, so that there is a clutching map
\begin{diagram}
\kappa:& \til \calm^{\gamma_1}_G\cross \til \calm^{\gamma_2}_G & \rto & \til \calm^\gamma_G.
\end{diagram}
\begin{alphabetize}
\item There is a canonical isomorphism of sheaves on $(\til
  \calm^{\gamma_1}_G \cross \til \calm^{\gamma_2}_G)\bc k$,
\begin{equation}
\label{eqclutchsum}
(\kappa\bc k)^* \pic^0(\calc^\gamma)[\ell] \iso
\pic^0(\calc^{\gamma_1})[\ell] \cross \pic^0(\calc^{\gamma_2})[\ell].
\end{equation}

\item Suppose $s_i \in \til \calm^{\gamma_i}_G(k)$ for $i= 1,
  2$, and  
  let $s = \kappa(s_1,s_2)$.  After base change to $k$, there is a commutative diagram:
\begin{diagram}
\mono_\ell(\til\calm^{\gamma_1}_G, s_1) \cross
\mono_\ell(\til\calm^{\gamma_2}_G, s_2) & \rto^\alpha_{\sim} &
\mono_\ell(\kappa(\til \calm^{\gamma_1}_G \cross \til
\calm^{\gamma_2}_G), s) & \rinject &
\mono_\ell(\til\calm^\gamma_G,s) \\
\dinject &&& \rdinject &\dinject\\
\aut(\pic^0(\calc^{\gamma_1})[\ell]_{s_1}) \cross
\aut(\pic^0(\calc^{\gamma_2})[\ell]_{s_2})
&& \rinject^\delta &&
\aut(\pic^0(\calc^{\gamma})[\ell]_{s})
\end{diagram}

\end{alphabetize}
\end{lemma}

\begin{proof}
Let $C \ra S$ be a stable curve.  Suppose that $C$
is the union of two (not necessarily irreducible) proper,
connected $S$-curves $C_1$ and $C_2$ which intersect along a unique section.  If $\call$ is a line bundle on
$C$ of degree zero, then for $i = 1, 2$ the restriction
$\call\rest{C_i}$ is a line bundle of degree zero on $C_i$ \cite[9.1.2
and 9.2.13]{blr}.   Thus, there is a morphism of group functors
$(\kappa\bc k)^*\pic^0(C^\gamma) \ra \pic^0(C^{\gamma_1}) \cross
\pic^0(C^{\gamma_2})$, which restricts to a morphism on
$\ell$-torsion.  This is an isomorphism on stalks \cite[9.2.8]{blr},
and thus an isomorphism of sheaves.

For part (b), the homomorphism $\delta$ is induced by the isomorphism
in Equation \eqref{eqclutchsum}.  
For $i = 1, 2$, let $C_i \ra S_i$ denote the cover $\calc^{\gamma_i} \bc k \ra \til \calm_G^{\gamma_i} \bc k$. 
To define $\alpha$, consider the kernel of the representation $\rho_{C_i \ra S_i,s_i, \ell}$.
Since it is an open normal subgroup of the fundamental group, this kernel 
defines an irreducible Galois \'etale
cover $Y_i \ra S_i$.
Let $H_i$ be its Galois group. 
Then $Y_1
  \cross_k Y_2$ is an irreducible $H_1\cross H_2$ cover of $S_1 \cross
  S_2$.
Using this and part (a), the cover 
$\pic^0(\calc^\gamma)[\ell]\rest{\kappa(S_1\cross S_2)}$ 
is trivialized by a Galois extension with group $H_1 \cross H_2$.
Therefore, $\alpha$ is an isomorphism. 
\end{proof}

\begin{lemma}
\label{lemmaxgp}
 Let $\ff$ be a finite field of odd characteristic, and let $V/\ff$ be
 a $g$-dimensional vector space.  Suppose that
\begin{equation}
\label{eqsumdecomp}
V = V_1 \oplus V_2 \oplus V_3,
\end{equation}
and that $g_i := \dim V_i$ is positive for $i=1,2,3$.  Let $H$ be
a subgroup of $\gl(V)$.

\begin{alphabetize}
\item\label{lemmaxsl} If $H\subseteq \sl(V)$ contains both $\sl(V_1\oplus V_2)$ and
  $\sl(V_2\oplus V_3)$, then $H = \sl(V)$. 

\item Suppose that $V$ is equipped with a nondegenerate pairing $\ang{\cdot,\cdot}$ and that for $i\not = j$, $\ang{V_i, V_j} = (0)$.
\begin{romanize}
\item\label{lemmaxsp} If the pairing is symplectic, and if $H\subseteq \sp(V)$ 
  contains both $\sp(V_1\oplus V_2)$ and $\sp(V_2\oplus V_3)$, then $H =
  \sp(V)$.

\item\label{lemmaxsu} If the pairing is Hermitian, and if $H\subseteq \su(V)$ contains both
  $\su(V_1\oplus V_2)$ and $\su(V_2\oplus V_3)$, then $H =\su(V)$.
\end{romanize}
\end{alphabetize}
\end{lemma}

\begin{proof} 
We use  \eqref{eqsumdecomp} to identify $\gl(V_i)$
with a subgroup of $\gl(V)$.

Our proof of (a) is Lie-theoretic.  Since $\sl(V)$ is split, the roots
and Weyl group of $\sl(V)$ are the same as those of
$\sl(V\tensor\bar\ff)$ \cite[1.18]{carterfglt}.  Moreover,
$\sl(V)$ has a split BN-pair, and we show $H = \sl(V)$ by successively
showing that $H$ contains a maximal torus $T$ of $\sl(V)$, the
associated Weyl group, and all root subgroups.

Choose coordinates $e_1, \ldots, e_g$ so that $V_1$ is the span of
$\st{e_1, \ldots, e_{g_1}}$; $V_2$ is the span of $\st{e_{g_1+1},
  \ldots, e_{g_1+g_2}}$; and $V_3$ is the span of $\st{e_{g_1+g_2+1},
  \ldots, e_g}$.  

First, $H$ contains the maximal split torus $T$ of $\sl(V)$
consisting of all diagonal matrices $\diag(\nu_1, \ldots, \nu_g)$ such that
$\prod_j \nu_j = 1$.  To see this, choose $j$ such that $e_j\in
V_2$. Given $\nu = \diag(\nu_1, \ldots, \nu_g) \in T$, we use the fact
that $\nu_j = (\prod_{i\not= j}\nu)\inv$ to write
\begin{equation*}
\nu =\diag(\nu_1, \ldots, \nu_{j-1}, (\prod_{i=1}^{j-1}\nu_i)\inv, 1,
\ldots, 1) \cdot \diag( 1, \ldots, 1,(\prod_{i=j+1}^g \nu_i)\inv,
\nu_{j+1}, \ldots, \nu_g).
\end{equation*}

Second, $H$ contains $N_{\sl(V)}(T)$, the normalizer of $T$ in
$\sl(V)$.  It suffices to show $N_H(T)/T = N_{\sl(V)}(T)/T$.  The
normalizer $N_{\sl(V)}(T)$ is the set of matrices of determinant one with exactly one
nonzero entry in each row and in each column, and the quotient
$N_{\sl(V)}(T)/T$ is isomorphic to $\sym(\st{1, \ldots, g})$.   Under
this identification, 
$$N_{\sl(V_1 \oplus V_2)}(T\cap \sl(V_1\oplus V_2))/(T\cap \sl(V_1\oplus V_2)) =
\sym(\st{1, \ldots, g_1+g_2});$$
$$N_{\sl(V_2\oplus V_3)}(T\cap \sl(V_2\oplus V_3))/(T\cap \sl(V_2\oplus V_3)) = \sym(\st{g_1+1,
  \ldots, g}).$$  
Since $H$ contains $\sl(V_1\oplus V_2)$ and
$\sl(V_2\oplus V_3)$, it contains $N_{\sl(V)}(T)$.

Third, let $\Delta$ be the canonical set of simple roots associated with
$\End(V) \iso \mat_n(V)$.  For each $ i \in \st{1,
  \ldots, g-1}$ there is a root $\alpha_i\in \Delta$.  The associated root group
is the unipotent group $U_{\alpha_i}$; the nontrivial elements of $U_{\alpha_i}$ are unipotent
matrices whose only nonzero offdiagonal entry is at $(i,i+1)$.  Clearly, $H$
contains each $U_{\alpha_i}$.

Finally, the Weyl group $N_{\sl(V)}(T)/T$, acting on the set of simple roots
$\Delta$, generates the set of all roots $\Phi$. Since $H$ contains
$N_{\sl(V)}(T)$, it therefore contains all root 
groups $U_\alpha$ where $\alpha\in \Phi$.  By the Bruhat 
decomposition, any element of $\sl(V)$ is a product of a
diagonal matrix and members of the various root groups.  Therefore, $H
= \sl(V)$.

For (b) part (i), after relabeling if necessary, we may assume that
$V_1\oplus V_2\not = V_3$.  Moreover, the hypothesis
ensures that $H$ contains $\sp(V_3)$, and thus it contains $\sp(V_1\oplus V_2)
\oplus \sp(V_3)$.  It is known that $\sp(V_1\oplus V_2) \oplus
\sp(V_3)$ is a maximal subgroup of $\sp(V)$ \cite[Thm.\ 3.2]{king81}. 
Since $\sp(V_2\oplus V_3)\not\subset
\sp(V_1\oplus V_2) \oplus \sp(V_3)$, the subgroup $H$ must equal $\sp(V)$.

The same argument proves (b) part (ii); as group-theoretic input, we use
the fact \cite[p.373]{king81} that $\su(V_1\oplus V_2)\oplus \su(V_3)$ is a
maximal subgroup of $\su(V)$.
\end{proof}

\subsection{Monodromy of hyperelliptic curves}

We show that the mod-$\ell$ monodromy of the tautological family of
hyperelliptic curves of genus $g$ is the full symplectic group
$\sp_{2g}(\integ/\ell)$.  Recall (Section \ref{subsechecomp}) that for
$g\in\nat$, the moduli space $\til\calh_g$ of labeled hyperelliptic
curves of genus $g$ is irreducible.  We will use Lemma
\ref{lemdisjointhe} to relate $\mono_\ell(\til\calh_g)$ to
$\mono_\ell(\bar\calh_g)$.

\begin{lemma}
\label{lemdisjointhe}
Let $\ell$ be an odd prime.  Let $S$ be an irreducible scheme with a
primitive $2\ell\th$ root of unity, and let $(X \ra S, \phi)$ be a principally
polarized abelian scheme.  Let $Y \ra S$ be a Galois \'etale cover.
Suppose that $\mono_\ell(X \ra S) \iso
\sp_{2g}(\integ/\ell)$ and that the groups $\gal(Y/S)$ and
$\sp_{2g}(\integ/\ell)$ have no common nontrivial quotients.  Then $\mono_\ell(
X \cross_S Y \ra Y) \iso \sp_{2g}(\integ/\ell)$.
\end{lemma}

\begin{proof}
We equip $(\integ/\ell)^{2g}$ with the standard symplectic pairing
$\ang{\cdot,\cdot}_\std$.  Let  
\begin{equation*}
I_\ell := \Isom( (X[\ell], \ang{\cdot,\cdot}_\phi),
((\integ/\ell)^{2g},\ang{\cdot,\cdot}_{\std})).
\end{equation*}
The hypothesis on $S$ implies that $I_\ell \ra S$ is an \'etale Galois cover with group
$\sp_{2g}(\integ/\ell)$.
The hypothesis on $\mono_\ell(X \ra S)$ implies that $I_\ell$ is irreducible.

To prove the lemma, we must show that $I_\ell \cross_SY$ is
irreducible.  Equivalently, we must show that $I_\ell \ra S$ and $Y
\ra S$ are disjoint.  Now, $I_\ell \cross_S Y \to S$ is a (possibly
reducible) \'etale Galois cover with group $\sp_{2g}(\integ/\ell)
\cross \gal(Y/S)$. If $Z \ra S$ is any common quotient of $I_\ell \ra
S$ and $Y \ra S$, then so is each conjugate $Z^\tau \to S$ for $\tau
\in \gal(I_\ell \cross_SY/S)$.  Therefore, the compositum $\til Z$ of
all such conjugates $Z^\tau$ is also a common quotient of $I_\ell \ra
S$ and $Y \ra S$.  It thus suffices to show that there is no
nontrivial Galois cover $\til Z \ra S$ which is a quotient of $I_\ell
\ra S$ and $Y \ra S$. This last claim is guaranteed by the
group-theoretic hypothesis.  \end{proof}

We compute the mod-$\ell$ monodromy of the moduli
space of hyperelliptic curves.

\begin{theorem}
\label{thhe}
Let $\ell$ be an odd prime, and let $k$ be an algebraically closed field in which $2\ell$ is invertible.
For $g\ge 1$, $\mono_\ell(\til\calh_g\bc k)\iso \mono_\ell(\calh_g\bc k) \simeq
\sp_{2g}(\integ/\ell)$.
\end{theorem}

\begin{proof}
Since $\calh_g\bc k$ is open and dense in $\bar\calh_g\bc k$, which is dominated
by $\til\calh_g\bc k$, it suffices to show that $\mono_\ell(\til\calh_g\bc
k) \iso \sp_{2g}(\integ/\ell)$.

Our proof is by induction on $g$.
If $g=1$ or $g= 2$, then every curve of genus $g$ is hyperelliptic,
so that $\bar\calm_g\bc k$ and $\bar\calh_g\bc k$ coincide.  Therefore,
$\mono_\ell(\bar\calh_g\bc k)\iso \mono_\ell(\bar\calm_g\bc k)$.   By
\cite[5.15-5.16]{delignemumford}, $\mono_\ell(\bar\calm_g\bc k) \iso
\sp_{2g}(\integ/\ell)$.

We conclude the base case $g=1$ and $g=2$ by applying Lemma \ref{lemdisjointhe} to the
cover $\til \calh_g \ra \bar \calh_g$.  This cover is \'etale and Galois, with Galois
group $\sym(2g+2)$.
This group and $\sp_{2g}(\integ/\ell)$ have no common nontrivial
quotient.  (To see this, recall that if $\ell$ is odd then the projective
symplectic group $\psp_{2g}(\integ/\ell)$ is simple, except that
$\psp_{2}(\integ/3) \iso A_4$.  Neither $\psp_{2g}(\integ/\ell)$ nor
any quotient of $A_4$ is a nontrivial quotient of $\sym(4)$ or $\sym(6)$.)
By Lemma \ref{lemdisjointhe}, $\mono_\ell(\til\calh_g\bc
k) \iso \sp_{2g}(\integ/\ell)$.

We now assume that $g \ge 3$ and that $\mono_\ell(\til\calh_{g'}\bc k) \iso
\sp_{2g'}(\integ/\ell)$ for $1 \le g' < g$.  By Lemma
\ref{lemdegenhe}, $\til\calh_g \bc k$
degenerates to $\Delta_{1,1}$, so that there is a
diagram \eqref{diagclutch}:
\begin{diagram}
(\til\calh_1 \cross \til \calh_{g-2} \cross \til \calh_1)\bc k &
\rto^{\kappa_R} &
(\til \calh_1 \cross \til \calh_{g-1})\bc k \\
\dto<{\kappa_L} & \rdto>\kappa & \dto \\
(\til \calh_{g-1} \cross \til \calh_1) \bc k& \rto & \til \calh_g \bc k.
\end{diagram}
Fix a base point $(s_1,s_2,s_3) \in (\til \calh_1 \cross \til
\calh_{g-2} \cross \til \calh_1)(k)$, and let $s =
\kappa(s_1,s_2,s_3)$.  For $h \in \NN$, let $C^{h} \to \til\calh_{h}$
be the tautological labeled $\ZZ/2$-curve of genus $h$.  Let $V =
\pic^0(\calc^g)[\ell]_s$, and for $i=1,2,3$ let $V_i =
\pic^0(\calc^{g_i})[\ell]_{s_i}$.  Each of these is a
$\integ/\ell$-vector space equipped with a symplectic form.  By Lemma
\ref{lemclutchmono}(a), there is an isomorphism of symplectic
$\integ/\ell$-vector spaces $V \iso V_1 \oplus V_2 \oplus V_3$.

Let $\til\calb_R$ (resp.\ $\til\calb_L$) be the image of $\kappa_R$ (resp.\ $\kappa_L$).    
Using the decomposition in Equation \eqref{eqclutchsum},
we have inclusions:
\begin{equation}
\label{eqhecontain}
\begin{split}
\mono_\ell(\til\calb_R\bc k) & \subseteq \sp(V_1) \oplus \sp(V_2 \oplus
V_3); \\
\mono_\ell(\til\calb_L\bc k) & \subseteq \sp(V_1\oplus V_2) \oplus \sp(V_3).
\end{split}
\end{equation}

By the induction hypothesis, the inclusions in \eqref{eqhecontain} are
equalities.  By Lemma \ref{lemclutchmono}(b), 
$\mono_\ell(\til \calh_g\bc k)$ contains $\sp(V_1)\oplus
  \sp(V_2\oplus V_3)$ and $\sp(V_1\oplus V_2) \oplus \sp(V_3)$.
By Lemma \ref{lemmaxgp}b(i), 
$\mono_\ell(\til\calh_g) \iso \sp_{2g}(\integ/\ell)$.
\end{proof}

\begin{corollary}
\label{corladiche}
Let $\ell$ be an odd prime, and let $k$ be an algebraically closed
field in which $2\ell$ is invertible.  For $g \ge 1$,
$\mono_{\integ_\ell}(\til \calh_g\bc k) \iso
\mono_{\integ_\ell}(\calh_g \bc k) \iso\sp_{2g}(\integ_\ell)$.
\end{corollary}

\begin{proof}
By Theorem \ref{thhe}, $\mono_\ell(\til\calh_g\bc k) \iso
\sp_{2g}(\integ/\ell)$.  By construction, there is a surjection
$\mono_{\integ_\ell}(\til\calh_g\bc k) \ra \mono_\ell(\til\calh_g \bc
k)$.  
Since the composition $\mono_{\integ_\ell}(\til \calh_g\bc k)
\inject \sp_{2g}(\integ_\ell) \ra \sp_{2g}(\integ/\ell)$ is
surjective, a standard group theory argument (e.g., \cite[1.3]{vasiusurj})
shows that $\mono_{\integ_\ell}(\til\calh_g \bc k) \iso
\sp_{2g}(\integ_\ell)$.
As in the proof of Theorem \ref{thhe}, $\mono_{\integ_\ell}(\calh_g \bc k) \iso \sp_{2g}(\integ_\ell)$
as well.
\end{proof}

\begin{corollary} \label{Chypabsirr}
Let $p \not = 2$ and let $g\in\nat$.  Then there exists a
smooth hyperelliptic curve of genus $g$ defined over $\bar\ff_p$ whose Jacobian is
absolutely simple.
\end{corollary}

\begin{proof}
Let $\ell$ be an odd prime distinct from $p$. By Corollary
\ref{corladiche},
\begin{equation*}
\sp_{2g}(\integ_\ell) \iso \mono_{\integ_\ell}(\calh_g\bc {\bar\ff_p})
\subseteq \mono_{\integ_\ell}(\calh_g\bc\ff_p) \subseteq
\gsp_{2g}(\integ_\ell).
\end{equation*}
Moreover,
$\mono_{\integ_\ell}(\calh_g\bc\ff_p)/\mono_{\integ_\ell}(\calh_g\bc{\bar\ff_p})
\iso \gal(\ff_{p,\ell^\infty}/\ff_p)$, which has finite index in
$\aut(\mmu_{\ell^\infty}(\bar\ff_p)) \iso
\gsp_{2g}(\integ_\ell)/\sp_{2g}(\integ_\ell)$.  Therefore, 
$\mono_{\integ_\ell}(\calh_g\bc \ff_p)$ is an open subgroup of
$\gsp_{2g}(\integ_\ell)$; by Borel's density theorem \cite[Th. 4.10]{platonovrapinchuk}, $\mono_{\rat_\ell}(\calh_g\bc \ff_p)
\iso \gsp_{2g}(\rat_\ell)$.  The claim now follows from
\cite[Prop. 4]{chaioort01}.
\end{proof}

\subsection{Monodromy of trielliptic curves}

We now compute the monodromy groups of tautological families
$\calc^\gamma \ra \til\calt^\gamma$ of labeled trielliptic curves. 
The Jacobians of trielliptic curves admit an action by
$\integ[\zeta_3]$.  This places a constraint on
$\mono_\ell(\til\calt^\gamma)$; in Theorem \ref{thtri}, we show that
this is the only constraint. 
We need more notation concerning $\integ[\zeta_3$]-actions and unitary groups
to describe the monodromy group precisely.

\subsubsection{Unitary groups}

Let $S$ be an irreducible scheme over $\spec\integ[1/3,\zeta_3]$, and
let $X \ra S$ be an abelian scheme of relative dimension $g$ equipped
with an action $\iota:\integ[\zeta_3] \ra \End_S(X)$.  Then $\Lie(X)$
is a locally free $\integ[\zeta_3]\tensor \calo_S$-module of some rank
$(r,s)$, where $r+s = g$.  We call $(r,s)$ the signature of the action
of $\integ[\zeta_3]$ on $X$.  If $(C/S,\iota_0)$ is a trielliptic curve,
then $\Lie(\pic^0(C))$ is the $\calo_S$-linear dual of
$\psi_*\Omega^1_{C/S}$ and the signature of $\pic^0(C)$ is the
same as that of $C$.

Let $V_{(r,s)}$ be a free $\integ[\zeta_3]$-module of rank
$g$, equipped with a $\integ[\zeta_3]$-linear pairing
$\ang{\cdot,\cdot}$ of signature $(r,s)$.  Let $\gu_{(r,s)}$ be the
$\integ[1/3, \zeta_3]$-group scheme of similitudes of the pair
$(V_{(r,s)}, \ang{\cdot,\cdot})$, and let $\g_{(r,s)}$ be the
restriction of scalars of $\gu_{(r,s)}$ to $\integ[1/3]$.  Let
$\su_{(r,s)} \subset \gu_{(r,s)}$ be the sub-group scheme of elements
of determinant one, and let $\sg_{(r,s)}$ be the restriction of
scalars of $\su_{(r,s)}$ to $\integ[1/3]$.    Concretely, for any
$\integ[1/3]$-algebra $R$, 
\begin{equation}
\label{eqdefsg}
\sg_{(r,s)}(R) = \st{ \tau \in \aut( V_{(r,s)}\tensor_{\integ[1/3]}R,
  \ang{\cdot,\cdot}): \det(\tau) = 1  }.
\end{equation}
In the abstract, we used $\su_{(r,s)}(\integ/\ell\tensor\integ[\zeta_3])$ to denote
$\sg_{(r,s)}(\integ/\ell)$.  The signature
condition implies that $\sg_{(r,s)}(\real)$ is isomorphic to the
complex special unitary group $\su(r,s)$.

The behavior of $\sg_{(r,s)}$ at finite 
primes $\ell \geq 5$ depends on their splitting in $\integ[\zeta_3]$.
Specifically, if $\ell$ is inert in $\integ[\zeta_3]$,
then $\sg_{(r,s)}(\integ/\ell) \iso \su_g(\ff_{\ell^2})$.
Alternatively, suppose that $\ell$ splits in $\integ[\zeta_3]$, and
let $\lambda$ and $\bar\lambda$ be the two primes of $\integ[\zeta_3]$
lying over $\ell$.  The factorization $\ell = \lambda \cdot
\bar\lambda$ yields a factorization $\integ[\zeta_3]\tensor\integ/\ell
\iso \integ/\ell \oplus \integ/\ell$.  This induces a decomposition
$V_{(r,s)}\tensor \integ/\ell \iso V_{(r,s)}(\lambda)\oplus V_{(r,s)}(\bar\lambda)$, where
$V_{(r,s)}(\lambda)$ and $V_{(r,s)}(\bar\lambda)$ are $g$-dimensional
$\integ/\ell$-vector spaces.  Moreover, the inner product $\ang{\cdot,\cdot}$ restricts
to a perfect pairing between $V_{(r,s)}(\lambda)$ and $V_{(r,s)}(\bar\lambda)$.  Then
$\sg_{(r,s)}(\integ/\ell)\subset \aut(V_{(r,s)}(\lambda)\cross
V_{(r,s)}(\bar\lambda))$ is the image of $\sl_g(\integ/\ell)$ embedded
as $\tau \mapsto \tau \cross (^t\tau)\inv$ (where $^t \tau$ denotes the transpose of $\tau$).    
In either case, the
isomorphism class of $\sg_{(r,s)}(\integ/\ell)$ depends only on $g$
and the splitting of $\ell$ in $\integ[\zeta_3]$, and not on the 
signature $(r,s)$.

\subsubsection{Calculation of trielliptic monodromy}

Let $\til\calt^\gamma$ be any component of the moduli space of
labeled trielliptic curves.  The signature of the action of $\integ[\zeta_3]$
on the Lie algebra of the relative Jacobian is locally constant, and we denote
it by $(r_\gamma, s_\gamma)$ (see Lemma \ref{Lsig}).  
The monodromy group $\mono_\ell(\til
\calt^\gamma\bc k)$ must preserve both the symplectic pairing and the
$\integ[\zeta_3]$-action on the $\ell$-torsion of
$\pic^0(\calc^\gamma\bc k)$.  This means that there is an
inclusion $\mono_\ell(\til \calt^\gamma\bc k) \subseteq
\sg_{(r_\gamma,s_\gamma)}(\integ/\ell)$.  We show that this is actually an
isomorphism.

As a basis for induction, we compare the moduli space of
trielliptic curves of genus three to the Picard modular variety, which
is a
component of a Shimura variety of PEL type.  Let $\calsh_{(2,1)}$ be the moduli stack
parametrizing data $(X/S, \iota, \phi)$ where $(X/S, \phi)$ is a principally polarized
abelian scheme of relative dimension $3$ and $\iota:\integ[\zeta_3]\ra
\End(X)$ 
satisfies signature and involution constraints.  Specifically, we
require that $\Lie(X)$ be a locally free $\integ[\zeta_3]\tensor
\calo_S$-module of signature $(2,1)$ and that $\iota$ take complex
conjugation on $\integ[\zeta_3]$ to the Rosati involution on
$\End(X)$.  Then $\calsh_{(2,1)}$ is the Shimura variety associated to
the reductive group $\g_{(2,1)}$.

We also
consider the moduli stack $\calsh_{(2,1),\ell}$ parametrizing data $(X/S, \iota,
\phi,\xi)$ where $(X/S,\iota,\phi)$ is as above and $\xi$ is a
principal level $\ell$ structure.   More precisely, $\xi$ is a
$\integ[\zeta_3]$-linear isomorphism  $\xi: X[\ell] \ra (V_{(2,1)}\tensor
\integ/\ell)\tensor \calo_S$ compatible with the given pairings.  The forgetful
functor $\calsh_{(2,1),\ell}  \ra \calsh_{(2,1)}$ induces a Galois
cover of stacks, with covering group $\sg_{(2,1)}(\integ/\ell)$
\cite[1.4]{gordon92}.
Note
that, as in \cite[IV.6.1]{faltingschai} but in contrast to \cite{bellaiche,larsenthesis}, we have implicitly chosen an
isomorphism $\mu_\ell \ra \integ/\ell$.     If
we had not chosen this isomorphism,
the covering group would be $\g_{(2,1)}(\integ/\ell)$; compare \cite[IV.6.12]{faltingschai}.

\begin{lemma}
\label{lemdisjointtri}
Let $\ell\ge 5$ be prime, and  let $S$ be an irreducible scheme with a primitive $3\ell\th$ root of unity.
Let $(X/S, \iota,\phi)\in \calsh_{(2,1)}(S)$ be a principally
polarized abelian scheme with $\integ[\zeta_3]$-action.  
Let $Y \ra S$ be an \'etale Galois cover.  Suppose that $\mono_\ell(X \ra S)
\iso \sg_{(2,1)}(\integ/\ell)$ and that the groups $\gal(Y/S)$ and
$\sg_{(2,1)}(\integ/\ell)$ have no common nontrivial quotients.  Then
$\mono_\ell(X \cross_S Y \ra Y) \iso \sg_{(2,1)}(\integ/\ell)$.
\end{lemma}

\begin{proof}
The proof is exactly the same as that of Lemma \ref{lemdisjointhe},
except that it involves 
$$I_{\integ[\zeta_3],\ell} := \Isom_{\integ[\zeta_3]}((X[\ell],
\ang{\cdot,\cdot}_\phi), (V_{2,1}\tensor\integ/\ell,
\ang{\cdot,\cdot})).$$  
Then $I_{\integ[\zeta_3],\ell} \to S$ is \'etale and Galois with group
$\su(V_{(2,1)}\tensor\integ/\ell) \simeq  \sg_{(2,1)}(\integ/\ell)$.  
The hypothesis on $\mono_\ell(X \ra S)$ is equivalent to the irreducibility of $I_{\integ[\zeta_3,],\ell}$.
\end{proof}

\begin{lemma}
\label{lemtribase}
Let $\ell\ge 5$ be prime, and let $k$ be an algebraically closed field
in which $3\ell$ is invertible.  
Let
$\til\calt^\gamma$ be any component of the moduli space of
trielliptic curves.  If $g(\gamma) = 3$, then 
$\mono_\ell(\til\calt^\gamma\bc k) \iso \sg_{(r_\gamma,s_\gamma)}(\integ/\ell)$.
\end{lemma}
\begin{proof}

We start by computing the monodromy group of the Picard modular
variety $\calsh_{(2,1)}$.
Both $\calsh_{(2,1)}\bc\cx$ and
$\calsh_{(2,1),\ell}\bc\cx$ are arithmetic quotients of the complex
ball, thus irreducible \cite[1.4]{gordon92}.  Moreover, the existence of smooth
arithmetic compactifications \cite[Section
  3]{larsenthesis} \cite[1.3.13]{bellaiche} and Zariski's
  connectedness theorem (see \cite[IV.5.10]{faltingschai}) imply that
  $\calsh_{(2,1)}\bc k$  and $\calsh_{(2,1),\ell}\bc k$ are also irreducible.  The irreducibility of the
  cover $\calsh_{(2,1),\ell}\bc k \ra \calsh_{(2,1)}\bc k$
  implies that the fundamental group of
  $\calsh_{(2,1)}\bc k$ acts transitively on the set of Hermitian
  $\integ[\zeta_3]\tensor\integ/\ell$-bases for the $\ell$-torsion of
  the tautological abelian scheme over $\calsh_{(2,1)}\bc k$.
Therefore, $\mono_\ell(\calsh_{(2,1)}\bc k) \iso
\sg_{(2,1)}(\integ/\ell)$.

We now consider $\til\calt^\gamma\bc k$.  Possibly after relabeling
\eqref{diagrelabelG}, by Lemma \ref{Lsig} we may assume that
$(r_\gamma,s_\gamma) = (2,1)$ and that the inertia type $\bar\gamma$
is $\st{1,1,1,1,2}$.    Moreover, since $\dim \calsh_{(2,1)}
\bc k= 2$, the Torelli map gives an inclusion of 
$\calt^{\bar\gamma}$ onto an open subset of $\calsh_{(2,1)}$.
Thus $\mono_\ell(\bar\calt^{\bar\gamma}\tensor k) \iso 
\mono_\ell(\calsh_{(2,1)}\tensor k)\iso \sg_{(2,1)}(\integ/\ell)$.

By Lemma \ref{lemforgetgalois}, $\til\calt^\gamma \ra
\bar\calt^{\bar\gamma}$ is Galois with group $\sym(4)$.
Since $\ell \ge 5$, the quotient of
$\sg_{(2,1)}(\integ/\ell)$ by its center is simple.  Thus,
$\sym(4)$ and $\sg_{(2,1)}(\integ/\ell)$  have no common nontrivial
quotient, and  $\mono_\ell(\til\calt^\gamma\bc k) \iso \sg_{(2,1)}(\integ/\ell)$
by Lemma \ref{lemdisjointtri}.
\end{proof}

\begin{theorem}
\label{thtri}
Let $\ell\ge 5$ be prime, and let $k$ be an algebraically closed
field in which $3\ell$ is invertible.
Let
$\til\calt^\gamma$ be any component of the moduli space of labeled
trielliptic curves.  If $g(\gamma) \ge 3$, then
$\mono_\ell(\til\calt^\gamma\bc k) \iso \mono_\ell(\calt^{\bar \gamma}\bc k) \iso \sg_{(r_\gamma,s_\gamma)}(\integ/\ell)$.
\end{theorem}

\begin{proof}
Since $\calt^{\bar \gamma}\bc k$ is open and dense in $\bar\calt^{\bar \gamma}\bc k$, which is dominated by $\til\calt^{\gamma}\bc k$, 
it is sufficient to show $\mono_\ell(\til\calt^\gamma\bc k) \iso \sg_{(r_\gamma,s_\gamma)}(\integ/\ell)$.
We proceed by induction on $g(\gamma)$.
The case $g(\gamma) = 3$ is supplied by Lemma \ref{lemtribase}.

Now let $\til\calt^\gamma$ be any component of $\til\calt$ with
$g = g(\gamma) \ge 4$, and
suppose the result is true for all components
$\til\calt^{\gamma'}$
with $3 \le g(\gamma') < g$.
By Proposition \ref{propdegentri}, $\til \calt^\gamma$ degenerates to $\Delta_{1,1}$.
This means that there are class vectors
$\gamma_1$, $\gamma_{2}$, $\gamma_3$, $\gamma_L$, and $\gamma_R$
such that $g(\gamma_1) = g(\gamma_3)
= 1$; $g(\gamma_2)=g-2$; and there is a diagram \eqref{diagclutch}:
\begin{diagram}
(\til \calt^{\gamma_1} \cross \til \calt^{\gamma_{2}} \cross \til
\calt^{\gamma_3})\bc k & \rto^{\kappa_R} & (\til \calt^{\gamma_1} \cross \til
\calt^{\gamma_R})\bc k \\
\dto<{\kappa_L} &\rdto^\kappa& \dto \\
(\til \calt^{\gamma_L} \cross \til \calt^{\gamma_3})\bc k & \rto &
\til \calt^\gamma\bc k.
\end{diagram}

Let $(s_1,s_2,s_3) \in (\til \calt^{\gamma_1} \cross \til \calt^{\gamma_{2}} \cross \til
\calt^{\gamma_3})(k)$, and let $s = \kappa(s_1,s_2,s_3)$.
Let $V =\pic^0(\calc^\gamma)[\ell]_s$, and for $i=1,2,3$ let $V_i=
\pic^0(\calc^{\gamma_i})[\ell]_{s_i}$.  Then $V \iso V_1 \oplus V_2 \oplus V_3$
is a decomposition of $V$ as a Hermitian $(\integ[\zeta_3]\tensor
k)$-module (Lemma \ref{lemclutchmono}(a)).  Let $\til \calb_R$ be the
image of $\kappa_R$, and let $\til \calb_L$ be the image of
$\kappa_L$.    Lemma \ref{lemclutchmono}(b) shows that $\mono_\ell(\til
\calt^\gamma\bc k, s)$ contains $\mono_\ell(\til \calb_L\bc k)$ and $\mono_\ell(\til \calb_R\bc k) $;
the inductive hypothesis shows
that these are $\su(V_1\oplus V_2)$ and $\su(V_2\oplus V_3)$,
respectively.

If $\ell$ is inert in $\integ[\zeta_3]$, Lemma
\ref{lemmaxgp}b(ii) implies that $\mono_\ell(\til\calt^\gamma\bc k) \iso
\su_g(\FF_{\ell^2}) \iso \sg_{(r_\gamma,s_\gamma)}(\integ/\ell)$.

Otherwise, if $\ell = \lambda \cdot \bar \lambda$ is split in
$\integ[\zeta_3]$, let $V(\lambda)$ be the eigenspace of $V$
corresponding to $\lambda$, and define $V_i(\lambda)$ analogously for
$i = 1,2,3$.
In this case, $\sg_{(r_\gamma,s_\gamma)}(\integ/\ell)$, $\su(V)$ and $\sl(V(\lambda))$ are isomorphic.   
By the inductive hypothesis,
the projection of $\mono_\ell(\til\calt^\gamma\bc k)$ to
$\sl(V(\lambda))$ contains $\sl(V_1(\lambda)\oplus V_2(\lambda))$ and
$\sl(V_2(\lambda)\oplus V_3(\lambda))$.  By Lemma \ref{lemmaxgp}(a),
we see that $\mono_\ell(\til\calt^\gamma\bc k) \iso \sl(V(\lambda)) \iso
\sg_{(r_\gamma,s_\gamma)}(\integ/\ell)$.
\end{proof}

\begin{corollary}
\label{corladictri}
Let $\ell$ be an odd prime, and let $k$ be an algebraically closed
field in which $3\ell$ is invertible.  Let $\til\calt^\gamma$ be any
component of the moduli space of trielliptic curves.  If $g(\gamma)
\ge 3$, then $\mono_{\integ_\ell}(\til\calt^\gamma\bc k) \iso
\mono_{\integ_\ell}(\calt^{\bar\gamma}\bc k) \iso
\sg_{(r_\gamma,s_\gamma)}(\integ_\ell)$.
\end{corollary}

\begin{proof}
The proof is parallel to that of Corollary \ref{corladiche}; any
subgroup of $\sg_{(r_\gamma,s_\gamma)}(\integ_\ell)$ which surjects
onto $\sg_{(r_\gamma,s_\gamma)}(\integ/\ell)$ is all of
$\sg_{(r_\gamma,s_\gamma)}(\integ_\ell)$.
\end{proof}

\begin{corollary} \label{Ctriabsirr}
Let $p \not = 3$.  Let $g \geq 3$ and let $(r,s)$ be a trielliptic signature for $g$ (Definition \ref{Dtrisig}).
Then there exists a smooth trielliptic curve defined over $\bar \ff_p$ with genus $g$ and signature $(r,s)$
whose Jacobian is absolutely simple.
\end{corollary}

\begin{proof}
  Let $\ff = \ff_{p^2}$.  Let $\ell$ be an odd prime distinct from $p$
  which is inert in $\integ[\zeta_3]$, and let $K_\ell =\rat_\ell(\zeta_3)$.
  By Lemma \ref{Lsig}, there is a
  class vector $\gamma$ whose inertia type has signature $(r,s)$.  As
  in the proof of Corollary \ref{Chypabsirr}, using Corollary
  \ref{corladictri} one sees that
  $\mono_{\rat_\ell}(\til\calt^\gamma\bc \ff)\iso\gu_{(r,s)}(K_\ell)$.

  Let $F$ be a CM field of degree $[F:\rat] = 2g$ which contains
  $\rat(\zeta_3)$ and is inert at $\ell$, and let $F_\ell = F
  \tensor\rat_\ell$.  There is a torus $H\subset 
  \gu_{(r,s)}(K_\ell)$ isomorphic to $F_\ell\units$.
  This torus is maximal since $[F:\rat] = 2g$.  The quotient of $H$ by
  the center of $\gu_{(r,s)}(K_\ell)$ is isogenous to the kernel
  $J$ of the norm map $F_\ell\units \ra \rat_\ell\units$.
  Since $F_\ell$ is a field, $J$ is anisotropic and $H$ is
  elliptic.  Finally, $H$ acts irreducibly on the Tate module of the
  Jacobian of the tautological curve $C^{\gamma}$, since
  $T_\ell(\pic^0(\calc^\gamma))\tensor\rat_\ell$ is a one-dimensional $F_\ell$-vector space.  The result then follows from
  \cite[Remark 5(i)]{chaioort01}.
\end{proof} 

\bibliographystyle{abbrv} 
\bibliography{jda}

\end{document}